\documentclass[10pt, leqno]{article}
\usepackage{amssymb}
\usepackage{amsmath}

\def\R{{\bf R}}
\def\F{{\bf F}}

\def\Pr{{\bf Pr}}

\def\inxx#1{}

\def\MYdef{\mathrel{\stackrel{\rm def}=}}

\newtheorem{theorem}{Theorem}[section]
\newtheorem{lemma}[theorem]{Lemma}
\newtheorem{corollary}[theorem]{Corollary}

\newtheorem{definition}[theorem]{Definition}

\numberwithin{equation}{section}

\newcommand\proof{{\bf Proof. }\nobreak\noindent}

\author{Anwar A. Irmatov}
\title{ On linear-combinatorial problems associated with subspaces spanned by $\{\pm 1\}$-vectors}

\date{}

\begin{document}

\renewcommand\refname{\centering  \sc References }

\maketitle
\begin{abstract}
A complete answer to the question about subspaces generated by $\{\pm 1\}$-vectors, which arose in the work of I.~Kanter and H.~Sompolinsky on associative memories, is given. More precisely, let vectors $v_1, \ldots , v_p,$ $p\leq n-1,$ be chosen at random uniformly and independently from $\{\pm 1\}^n \subset {\bf R}^n.$ Then the probability ${\mathbb P}(p, n)$ that 
$$span \ \langle v_1, \ldots , v_p \rangle \cap \left\{ \{\pm 1\}^n \setminus \{\pm v_1, \ldots , \pm v_p\}\right\} \ne \emptyset \ $$
\noindent is shown to be
$$4{p \choose 3}\left(\frac{3}{4}\right)^n + O\left(\left(\frac{5}{8} + o_n(1)\right)^n\right) \quad \mbox{as} \quad n\to \infty,$$
\noindent where the constant implied by the $O$-notation does not depend on $p$. The main term in this estimate is the probability that some 3 vectors $v_{j_1}, v_{j_2}, v_{j_3}$ of $v_j$, $j= 1, \ldots , p,$ have a linear combination that is a $\{\pm 1\}$-vector different from $\pm v_{j_1}, \pm v_{j_2}, \pm v_{j_3}. $

\bigskip


\end{abstract}

\let\thefootnote\relax


\section{Introduction.}

In the paper {\rm \cite{Odl}}, A.M.Odlyzko  gave a partial answer to the following question that arose in the paper {\rm \cite{KaSo}} on associative memories. Let vectors $v_1, \ldots , v_p$ be chosen at random uniformly and independently from 
$\{\pm 1\}^n \subset {\bf R}^n.$ What is the probability ${\mathbb P}(p, n)$ that the subspace spanned by $v_1, \ldots , v_p$ over reals contains a $\{\pm 1\}$-vector different from $\pm v_1, \ldots , \pm v_p$, i.e.
\begin{equation*}
span \ \langle v_1, \ldots , v_p \rangle \cap \left\{ \{\pm 1\}^n \setminus \{\pm v_1, \ldots , \pm v_p\}\right\} \ne \emptyset \  ?
\end{equation*}
G.Kalai and N.Linial conjectured (see {\rm \cite{Odl}}) that this probability is dominated by the probability ${\mathbb P}_3(p, n)$ that some 3 vectors $v_{j_1}, v_{j_2}, v_{j_3}$ of $v_j$, $j= 1, \ldots , p,$ have a linear combination that is a $\{\pm 1\}$-vector different from $\pm v_{j_1}, \pm v_{j_2}, \pm v_{j_3}, $ where 
\begin{equation}\label{eqP1}
{\mathbb P}_3(p, n)=  4{p \choose 3}\left(\frac{3}{4}\right)^n + O\left(p^4\left(\frac{5}{8}\right)^n\right) \quad \mbox{as} \quad n\to \infty,
\end{equation}
\noindent for $3\leq p \leq n.$
In the paper {\rm \cite{Odl}} this conjecture was proven for 
\begin{equation}\label{eqP2}
p \leq n - \frac{10n}{\ln n}.
\end{equation}
\bigskip 
\begin{theorem}\label{Pthm1}
(A.M.Odlyzko {\rm \cite{Odl}}) 
If $p \leq n - \frac{10n}{\ln n}$ and vectors $v_1, \ldots , v_p$ are chosen at random uniformly and independently from $\{\pm 1\}^n \subset {\bf R}^n$, then the probability ${\mathbb P}(p, n)$ that the subspace spanned by $v_1, \ldots , v_p$ over reals contains a $\{\pm 1\}$-vector different from $\pm v_1, \ldots , \pm v_p$ equals
\begin{equation}\label{eqP3}
{\mathbb P}(p, n)=  {\mathbb P}_3(p, n)+ O\left(\left(\frac{7}{10}\right)^n\right) \quad \mbox{as} \quad n\to \infty.
\end{equation}
The constants implied by the $O$-notation in (\ref{eqP1}) and (\ref{eqP3}) are independent of $p.$ 
\end{theorem}
\bigskip

The paper {\rm \cite{Odl}}  made a significant contribution to estimating the number of threshold functions $P(2, n).$ Namely, in the paper {\rm \cite{Zue}}, as a corollary of Theorem 1, T.Zaslavsky's formula {\rm \cite{Zas}} and G.-C.Rota's theorem on the inequality of the
 M\"obius function of a geometric lattice to zero (Theorem 4 {\rm \cite{Rot}}, p.357), the asymptotics of the logarithm of $P(2,n)$ was obtained.
In {\rm \cite{Ir1}}, the Theorem 1  was used together with the original $(A,B,C)$-construction to improve the lower bound of $P(2, n)$ obtained in {\rm \cite{Zue}}  by a factor of $\thicksim P(2, \left\lfloor \frac{7n\ln 2}{\ln n}\right\rfloor)$. In {\rm \cite{IrK}}, the Theorem 1 was generalized to the case of $E_K$-vectors, where 
$E_K = \{0, \pm 1, \ldots , \pm Q\}$ if $K=2Q+1$, and $E_K = \{\pm 1, \pm 3, \ldots , \pm (2Q-1)\}$ if $K=2Q,$  to obtain the asymptotics of the logarithm of the number of threshold functions of K-valued logic.

In {\rm \cite{Odl}}, A.M.Odlyzko noted that removing constraint (\ref{eqP2}) is an open and hard problem. The hardness of this problem is manifested in the fact (see {\rm \cite{KKS}}, p. 238), that if ${\mathbb P}(n-1, n)$ tends to zero as $ n\to \infty ,$ then one could obtain the asymptotics of $ P(2, n):$
\begin{equation}\label{eqP4}
P(2, n) \thicksim 2 {2^n-1 \choose n},\quad n\to \infty.
\end{equation}
We will show this fact (see Corollary \ref{corP1}) using the properties of $\eta^{\bigstar}_n$ function from the paper {\rm \cite{Ir4}}.

Also, in {\rm \cite{TV2}}, one can see an implicit connection between estimates of the probability ${\mathbb P}(n-1, n)$ and the probability $\mathbb{P}_n$ that a random Bernoulli $n\times n $  $\{\pm 1\}$-matrix $M_n$ is singular:
\begin{equation*}\label{eqPr}
\mathbb{P}_n \MYdef \Pr ( \det M_n = 0).
\end{equation*}
In this work, T.~Tao and V.~Vu investigate the properties of the combinatorial Grassmannian $Gr$, consisting of hyperplanes $V$ in an n-dimensional space over a finite field $\mathbb{F},$ $|\mathbb{F}|=p> n^{\frac{n}{2}}$, where $p$ is a prime number, such that $V= span \ \langle V  \cap  \{\pm 1\}^n \rangle$, and estimate $\mathbb{P}_n$ based on the formula 
\begin{equation*}
\mathbb{P}_n = 2^{o(n)}\sum_{V\in Gr,\  V \  is\ a\ non-trival\ hyperplane\ in\  \mathbb{F}^n}\mathbb{P}(A_V).
\end{equation*}
Here $A_V$ denotes the event that vectors $v_1, \ldots, v_n ,$ chosen at random uniformly and independently from $\{\pm 1\}^n ,$ span $V$, and $\mathbb{P}(A_V)$ denotes the probability of this event.  
We can essentially improve the upper bound of $\mathbb{P}(A_V)$ in the Small combinatorial dimension estimate Lemma from {\rm \cite{KKS}}  (see inequality (6) in Lemma 2.3. from {\rm \cite{TV2}}) if we substitute the probability $\mathbb{P}(v_1, \ldots , v_{n-1} \ \mbox{span} \ V)$ by the probability 
$$\mathbb{P}(\{v_1, \ldots , v_{n-1} \ \mbox{span} \ V\} \ \wedge \ \{V \cap \left\{ \{\pm 1\}^n \setminus \{\pm v_1, \ldots , \pm v_{n-1}\}\right\} \ne \emptyset\}).$$

On the other hand, J.~Kahn, J.~Koml\' os, and E.~Szemer\' edi in {\rm \cite{KKS}}, relying on a technique developed to estimate $\mathbb{P}_n,$ showed (corollary 4, {\rm \cite{KKS}}) that there is a constant $C$ such that the Theorem 1 is true for $p\leq n-C.$

In this paper, using the asymptotics of the probability $\mathbb{P}_n$ (Theorem 6, {\rm \cite{Ir4}})
\begin{equation}\label{eqP5}
\mathbb{P}_n \thicksim (n-1)^22^{1-n}, \quad n\to \infty
\end{equation}
and the lower bound for $P(2,n)$ (inequality (139) of the Theorem 7 {\rm \cite{Ir4}})
\begin{equation}\label{eqP6}
P(2, n) \geq  2\left[ 1- \frac{n^2}{2^n}\left( 1 + o\left(\frac{n^3}{2^n}\right) \right) \right]\dbinom{2^n-1}{n},
\end{equation}
we remove in the Theorem 1 the restriction (\ref{eqP2}) and prove the conjecture of G.Kalai, N.Linial, and A.M.Odlysko for $p\leq n-1.$

\bigskip 
\begin{theorem}\label{Pthm2}
{\bf (Main theorem)}
If $p \leq n - 1$ and vectors $v_1, \ldots , v_p$ are chosen at random uniformly and independently from $\{\pm 1\}^n \subset {\bf R}^n$, then the probability ${\mathbb P}(p, n)$ that the subspace spanned by $v_1, \ldots , v_p$ over reals contains a $\{\pm 1\}$-vector different from $\pm v_1, \ldots , \pm v_p$ equals
\begin{equation}\label{eqP7}
{\mathbb P}(p, n)=  {\mathbb P}_3(p, n)+ O\left(\left(\frac{5}{8}+ o_n(1)\right)^n\right) \quad \mbox{as} \quad n\to \infty.
\end{equation}
The constants implied by the $O$-notation in (\ref{eqP1}) and (\ref{eqP7}) are independent of $p.$ 
\end{theorem}
\bigskip

\section{Interdependence between ${\bf \eta}^{\bigstar}_n$-function, ${\mathbb P}(n, n+1)$, and the fraction of singular $(n+1)\times (n+1)$-$\{\pm 1\}$-matrices with different rows.}

Let $\langle H \rangle = \{\langle w_1\rangle, \ldots , \langle w_T\rangle \} \subset {\bf RP}^n$  be a subset of the $n$-dimensional projective space, where points $\langle w_i\rangle,$ $i=1, \ldots , T,$ are represented by lines $tw_i \subset \R^{n+1}$ and $w_i \in \R^{n+1},$ $i= 1, \ldots, T.$  Let $K^H$ be a simplicial compex defined as follows. The set of vertices of $K^H$ coincides with the set $\langle H \rangle$. A subset $\{\langle w_{i_1} \rangle , \ldots , \langle w_{i_s} \rangle \}$ of $\langle H \rangle$
forms a simplex of $K^H$ iff
\begin{equation*}
span \ \langle w_{i_1}, \ldots ,w_{i_s}\rangle \ne span \ \langle w_1, \ldots ,w_T\rangle .
\end{equation*}
We define on the set $2^{{\bf RP}^n}_{fin}$ of finite subsets of ${\bf RP}^n$ the function $\eta^{\bigstar}_n : 2^{{\bf RP}^n}_{fin} \to \mathbb{Z}_{\geq 0}$ by the following formula (see {\rm \cite{Ir4}}):
\begin{equation}\label{eqP8}
\eta^{\bigstar}_n(\langle H \rangle) = rank\ \tilde H_{n-1}(K^H;{\bf F}), \quad \langle H \rangle \subset {\bf RP}^n,
\end{equation}
\noindent if
\begin{equation*}
span \ \langle w_1, \ldots ,w_T\rangle = {\bf R}^{n+1},
\end{equation*}
\noindent and 
\begin{equation}\label{eqP9}
\eta^{\bigstar}_n(\langle H \rangle) = 0,
\end{equation}
\noindent if 
\begin{equation*}
span \ \langle w_1, \ldots ,w_T\rangle \ne {\bf R}^{n+1}.
\end{equation*}
\noindent Here $\tilde H_{n-1}(K^H;{\bf F})$ denotes the reduced homology group of the complex $K^{H}$ with coefficients in an arbitrary field $\F.$

Let us denote by $\langle H \rangle^{\times s}$, $s=1,\ldots, T$, the set of ordered collections $(\langle w_{i_1}\rangle , \ldots ,\langle w_{i_s} \rangle)$ of different $s$  elements from $\langle H \rangle$ and $\langle H \rangle^{\times s}_{\ne 0}\subset \langle H \rangle^{\times s},$ $\langle H \rangle^{\times s}_{= 0}\subset \langle H \rangle^{\times s}$ be the subsets 
\begin{equation}\label{eqP10}
\langle H\rangle ^{\times s}_{\ne 0} \MYdef \{(\langle w_{i_1}\rangle, \ldots ,\langle w_{i_s}\rangle) \in \langle H \rangle^{\times s} \mid dim \, span \ \langle w_{i_1}, \ldots , w_{i_s} \rangle  = s \}.
\end{equation}
\begin{equation}\label{eqP11}
\langle H\rangle ^{\times s}_{= 0} \MYdef \{(\langle w_{i_1}\rangle, \ldots ,\langle w_{i_s}\rangle) \in \langle H \rangle^{\times s} \mid dim \, span \ \langle w_{i_1}, \ldots , w_{i_s} \rangle  < s \}.
\end{equation}
\\
\noindent For any $W= (\langle w_{i_1} \rangle, \ldots , \langle w_{i_n}\rangle) \in \langle H \rangle^{\times n}$  and $ l = 1, \ldots, n,$ let 
\begin{equation}\label{eqP12}
\begin{split}
&L_l(W) \MYdef  span \   \langle w_{i_{n-l+1}}, \ldots , w_{i_n} \rangle \subset  {\bf R}^{n+1}; \\
&q^W_l(H) \MYdef |L_l(W) \cap \langle H \rangle|.
\end{split}
\end{equation}

\begin{definition}\label{Pdf1}
For any $W \in \langle H\rangle^{\times n},$ the ordered set of numbers
\begin{equation}\label{eqP13}
W(\langle H \rangle) \MYdef (q^W_n(H), q^W_{n-1}(H), \ldots, q^W_1(H))
\end{equation}
is called a combinatorial flag on $\langle H \rangle \subset {\bf RP}^n$ of the ordered set $W$.
\end{definition}
For the sake of simplicity, we will use the following notation:
\begin{equation}\label{eqP14}
W[H] \MYdef q^W_n(H)\cdot q^W_{n-1}(H) \cdots q^W_1(H).
\end{equation}
In {\rm \cite{Ir4}}, it was proven the following theorem.
\begin{theorem}\label{Pthm3}
{\bf ({\rm \cite{Ir4}}, {\rm \cite{Ir5}})}    For any $p=(p_1, \ldots , p_T),$ $p_i \in {\bf R}$, $i=1,\ldots, T$, such that $ \sum_{i=1}^T{p_i}=1,$ and subset $\langle H \rangle = \{\langle w_1\rangle, \ldots , \langle w_T\rangle \} \subset {\bf RP}^n,$ such that $span \ \langle H \rangle = {\bf R}^{n+1},$ the following equality is true: 
\begin{equation}\label{eqP15}
\eta^{\bigstar}_n(\langle H \rangle ) = \sum_{W\in \langle H\rangle^{\times n}_{\ne 0}}{\frac{1- p_{i_1} -p_{i_2} - \cdots - p_{i_{q_n^W}}}{W[H]}}.
\end{equation}
Here, the indices used in the numerator correspond to elements from 
$$L_n(W)\cap \langle H \rangle = \left\{\langle w_{i_1} \rangle , \ldots, \langle w_{i_n} \rangle, \ldots , \langle w_{i_{q_n^W}}\rangle \right\}.$$
\end{theorem}

Let 
\begin{equation}\label{eqP16}
E_n=\{ (1, b_1, \ldots , b_n) \mid \; b_i \in \{\pm 1 \} ,\; i=1, \ldots,  n \}
\end{equation}
\noindent and
\begin{equation}\label{eqP17}
\{E_n\}^{p}= \underbrace{E_n\times \cdots \times E_n}_{p}.
\end{equation}

We say that an ordered collection $W=(w_1, \ldots , w_p),$ $w_i \in \{\pm 1\}^{n+1},$ $i=1, \ldots , p,$ satisfies to KSO-condition, and we write $W\in KSO(p, n+1),$ iff 
\begin{equation}\label{eqP18}
span \ \langle w_1, \ldots , w_p \rangle \cap \left\{ \{\pm 1\}^{n+1} \setminus \{\pm w_1, \ldots , \pm w_p\}\right\} \ne \emptyset .
\end{equation}
\noindent Then
\begin{equation*}
{\mathbb P}(n, n+1)=\frac{|KSO(n, n+1)|}{2^{n(n+1)}}= 
\frac{|\{W\in \{E_n\}^n \mid W\in KSO(n, n+1)\}|}{2^{n^2}} ,
\end{equation*}
\noindent and
\begin{equation}\label{eqP19}
\begin{split}
&| KSO(n, n+1) \cap \langle E_n \rangle^{\times n}_{\ne 0}| =\\
&=|\{W\in \langle E_n\rangle^{\times n }_{\ne 0} \mid W\in KSO(n, n+1)\}| < 2^{n^2} {\mathbb P}(n, n+1).
\end{split}
\end{equation}

 \bigskip
\begin{corollary}\label{corP1}
If \ $\lim_{n\to \infty} {\mathbb P}(n, n+1) =0,$ then $P(2, n) \thicksim 2 {2^n-1 \choose n}, \  n\to \infty .$
\end{corollary}
\proof Let us apply Theorem \ref{Pthm3} to the set $H=E_n \subset \R^{n+1}$ and the collection of weights $p=(1, 0, \ldots, 0),$ where $w_1=(1, \ldots , 1) \in \R^{n+1}.$
Then
\begin{equation}\label{eqP20}
\eta^{\bigstar}_n(\langle E_n \rangle) = \Sigma_1  +\Sigma_2 - \Sigma_3 -\Sigma_4,
\end{equation}
\noindent where
\begin{equation*}
\Sigma_1 = \sum_{W\in \langle E_n \rangle^{\times n}_{\ne 0}, \ W\notin KSO(n, n+1)} \frac{1}{W[E_n]},
\end{equation*}
\begin{equation*}
\Sigma_2 = \sum_{W\in \langle E_n \rangle^{\times n}_{\ne 0}, \ W\in KSO(n, n+1)} \frac{1}{W[E_n]},
\end{equation*}
\begin{equation*}
\Sigma_3 = \sum_{W\in \langle E_n \rangle^{\times n}_{\ne 0}, \  w_1\notin \{W\},  \ w_1\in span \langle W \rangle} \frac{1}{W[E_n]},
\end{equation*}
\begin{equation*}
\Sigma_4 = \sum_{W\in \langle E_n \rangle^{\times n}_{\ne 0}, \  w_1\in \{W\}} \frac{1}{W[E_n]}.
\end{equation*}
\noindent From (\ref{eqP19}), we have
\begin{equation}\label{eqP21}
\begin{split}
\Sigma_1 &= \frac{|\langle E_n \rangle^{\times n}_{\ne 0}| - | KSO(n, n+1) \cap \langle E_n \rangle^{\times n}_{\ne 0}|}{n!} >\\
&>\frac{2^n\cdots (2^n-n+1) - 2^{n^2}{\mathbb P}_n - | KSO(n, n+1) \cap \langle E_n \rangle^{\times n}_{\ne 0}|}{n!} > \\
&>{2^n \choose n} - {2^n \choose n}\frac{2^{n^2}}{2^n(2^n-1)\cdots (2^n-n+1)}({\mathbb P}_n  +{\mathbb P}(n, n+1)),
\end{split}
\end{equation}

\begin{equation}\label{eqP22}
\Sigma_3 < {\mathbb P}_n(\{0,1\})\frac{2^{n^2}}{2^n(2^n-1)\cdots (2^n-n+1)} {2^n \choose n},
\end{equation}

\begin{equation}\label{eqP23}
\Sigma_4 < n\frac{(2^n-1)(2^n-2)\cdots (2^n-n+1)}{n!} = \dbinom{2^n-1}{n-1}.
\end{equation}
\noindent Here ${\mathbb P}_n(\{0,1\})$ denotes the probability that  a $(0, 1)$-$n\times n$-matrix, with entries chosen at random, uniformly, and independenly from $\{0,1\}$, is singular.
It follows from (\ref{eqP20}), (\ref{eqP21}), (\ref{eqP22}), and (\ref{eqP23}) that
\begin{equation}\label{eqP24}
\begin{split}
\eta^{\bigstar}_n(\langle E_n \rangle) >{2^n -1 \choose n}\left(1 - c_n)\right),
\end{split}
\end{equation}
\noindent where
\begin{equation*}
c_n = \frac{2^{n^2}}{(2^n-1)\cdots (2^n-n)}({\mathbb P}_n  +{\mathbb P}(n, n+1) +{\mathbb P}_n(\{0,1\})).
\end{equation*}
Taking into account the inequality (see the formulas (19) and (25) of {\rm \cite{Ir4}})
\begin{equation*}
P(2, n) \geq 2\eta^{\bigstar}_n(\langle E_n \rangle) ,
\end{equation*}
L.~Schl\"afli's upper bound (see the formula 2 in {\rm \cite{Ir4}} and {\rm \cite{Sch}})
\begin{equation*}
P(2, n) \leq 2 \sum_{i=0}^n {2^n-1 \choose i} ,
\end{equation*}
the given fact that ${\mathbb P}(n, n+1) \to 0,$ and the well known results ${\mathbb P}_n(\{0,1\}) \to 0,$ and ${\mathbb P}_n \to 0$ as $n\to \infty$ (see {\rm \cite{Ko1}}, {\rm \cite{Ko2}}, 
{\rm \cite{KKS}}, {\rm \cite{TV2}}, {\rm \cite{BVW}}, {\rm \cite{Tik}}, {\rm \cite{Ir4}}), we can conclude that 
\begin{equation*}
c_n \to 0 \ \  as \ \  n\to \infty ,
\end{equation*}
\noindent and
\begin{equation*}
P(2, n) \thicksim 2 {2^n-1 \choose n}, \  n\to \infty . 
\end{equation*}
\begin{flushright} {\sc Q.E.D.} \end{flushright}

We define $\delta_{n, k}$, $k= 1, \ldots, n+1,$ as
\begin{equation*}
\delta_{n, k} \MYdef \frac{|\langle E_n \rangle^{\times k}_{= 0}|}{|\langle E_n \rangle^{\times k}|}, \quad k=1, \ldots, n+1.
\end{equation*}

For $W = (\langle w_{i_1} \rangle, \ldots , \langle w_{i_n}\rangle) \in \langle E_{n}\rangle^{\times n}_{\ne 0},$ we use the following notations:
\begin{equation*}
\begin{split}
&L_n(W) \MYdef  span \   \langle w_{i_1}, \ldots , w_{i_n} \rangle \subset  {\bf R}^{n+1}; \\
&q_n^{W} \MYdef |L_n(W) \cap E_{n}|; \\
&E^m_{n} \MYdef  \left\{ W \in \langle E_{n} \rangle^{\times n}_{\ne 0} \mid q_n^{W} = n+m \right\}, \ \ m = 0, 1, \ldots , 2^{n-1} - n.
\end{split}
\end{equation*}
\bigskip

\begin{theorem}\label{Pthm4}
For sufficiently large $n,$ we have
\begin{equation}\label{eqP25}
\frac{|KSO(n, n+1) \cap \langle E_n \rangle^{\times n}_{\ne 0}|}{2^{n^2}} \leq \frac{n^2}{2^{n-1}} .
\end{equation}
\end{theorem}
\proof  Let us take a vector $w\in \R^{n+1}$ in general position to the set $E_n ,$ i.e. for any vectors $w_{i_1}, \ldots, w_{i_{n}} \in E_n \subset \R^{n+1},$  
$$w \notin span \langle w_{i_1}, \ldots, w_{i_{n}} \rangle .$$
 From the Theorem \ref{Pthm3} applyed to the set $H=\langle E_n\rangle \cup {\langle w \rangle } \subset {\bf RP}^{n}$ and the collection of weights $p(w) =1$ and $p(w_i) =0$ for $w_i \in E_n,$ $i=1, \ldots 2^n,$ we get
\begin{equation*}
\begin{split}
&\eta^{\bigstar}_n(\langle E_n\rangle \cup {\langle w \rangle })=\sum_{W\in \langle E_n \rangle^{\times n}_{\ne 0}} \frac{1}{W[E_n]} = \\
&=\sum_{m=0}^{2^{n-1}-n}\sum_{W\in E_n^m}\frac{1}{W[E_n]} = \frac{1}{n!}\sum_{m=0}^{2^{n-1}-n}|E_n^m| - \sum_{m=1}^{2^{n-1}-n}\sum_{W\in E_n^m}\left( \frac{1}{n!} - \frac{1}{W[E_n]}\right) \leq \\
&\left( \frac{1}{n!} - \frac{1}{W[E_n]} \geq \frac{1}{n!} -\frac{1}{(n-1)!(n+m)} = \frac{m}{n!(n+m)}\right)\\
&\leq \frac{1}{n!}\left(|\langle E_n \rangle^{\times n}|(1-\delta_{n,n}) -\frac{1}{n+1}\left|\cup_{m=1}^{2^{n-1}-n} \ E_n^m\right|\right),
\end{split}
\end{equation*}
or
\begin{equation}\label{eqP26}
\begin{split}
\eta^{\bigstar}_n(\langle E_n\rangle \cup {\langle w \rangle }) &\leq \dbinom{2^n}{n}(1-\delta_{n,n}) - \\
&-\frac{1}{n!(n+1)}\left|KSO(n, n+1)\cap \langle E_n \rangle^{\times n}_{\ne 0}\right|.
\end{split}
\end{equation}
From the Theorem 5 of the paper {\rm \cite{Ir4}}, we have (see the formula (135) in {\rm \cite{Ir4}}):
\begin{equation}\label{eqP27}
\eta^{\bigstar}_n(\langle E_{n} \rangle \cup \left\{\langle w \rangle \right\}) \geq \dbinom{2^n}{n}\left[ 1- \delta_{n, n} - \frac{n-1}{2^{n-1}}\left( 1 + o\left(\frac{n^3}{2^n}\right) \right) \right].
\end{equation}
Combining inequalities (\ref{eqP26}) and (\ref{eqP27}), we get 
\begin{equation}\label{eqP28}
 \left|KSO(n, n+1)\cap \langle E_n \rangle^{\times n}_{\ne 0}\right| \leq 2^n\cdots (2^n-n+1) \frac{n^2-1}{2^{n-1}}\left( 1 + o\left(\frac{n^3}{2^n}\right) \right).
\end{equation}
The Theorem \ref{Pthm4} follows from the inequality (\ref{eqP28}).
\begin{flushright} {\sc Q.E.D.} \end{flushright}

\bigskip

\section{Proof of the Main Theorem.}

\bigskip

Let ${\mathbb P}_m(p, n),$ $m\leq p \leq n-1,$ denote the probability that in the set of $p$ vectors $v_1, \ldots , v_p \in \{\pm 1 \}^n \subset \R^n,$ chosen at random uniformly and independently, there are some $m$ vectors $v_{j_1}, \ldots, v_{j_m}$ such that 
$$\alpha_1 v_{j_1} + \cdots + \alpha_m v_{j_m} \in \{\pm 1\}^n \ \ \mbox{for some} \ \ \alpha_1, \ldots, \alpha_m \in {\mathbb R}\setminus \{0\}.$$ 
Let $\mathcal{M}_m(p,n)$ denote the set of $(p\times n)$-$\{\pm 1\}$-matrices $M$ with linear independent rows $w_1, \ldots , w_p \in \{\pm 1\}^n$ satisfying the following property. There are a subset of $m$ rows  $w_{i_1}, \ldots w_{i_m}$ and some nonzero coefficients  
$ \alpha_1, \ldots, \alpha_m \in {\mathbb R}\setminus \{0\}$ such that 
$$\alpha_1 w_{i_1} + \cdots + \alpha_m w_{i_m} \in \{\pm 1\}^n.$$
Let $\mathcal{Q}(p,n)$ be the the set of $(p\times n)$-$\{\pm 1\}$-matrices $M$ with rank less than $p$ $(<p).$ 
Denote by $R_m(p, n)$ the probability that  a $(p\times n)$-$\{\pm 1\}$-matrix $M$ chosen at random belongs to $\mathcal{M}_m(p,n)$ and by ${\mathbb P}_{p,n}$ the probability that a $(p\times n)$-$\{\pm 1\}$-matrix $M$ chosen at random
has rank less than $p$ $(<p).$
Then 
\begin{equation*}
KSO(p, n) \subset \bigcup_{m=3}^p \mathcal{M}_m(p,n) \bigsqcup \mathcal{Q}(p,n),
\end{equation*}
\begin{equation}\label{eqP29}
R_m(p, n) \leq \dbinom{p}{m} R_m(m,n),
\end{equation}
\begin{equation}\label{eqP30}
{\mathbb P}_m(p, n) \leq R_m(p, n) + {\mathbb P}_{p, n} \leq \dbinom{p}{m} R_m(m,n) +{\mathbb P}_{p, n},
\end{equation}
and
\begin{equation}\label{eqP31}
{\mathbb P}(p, n) \leq \sum_{m=3}^p R_m(p, n) + {\mathbb P}_{p, n} \leq \sum_{m=3}^p\dbinom{p}{m} R_m(m,n) +{\mathbb P}_{p, n},
\end{equation}
It follows from {\rm \cite{Ir4}} (see Lemma 5 and the proof of the Theorem 6) that for $p= 1, \ldots , n$
\begin{equation}\label{eqP32}
{\mathbb P}_{p, n} \leq \frac{(p-1)^2}{2^{n-1}}(1 +o_n(1)).
\end{equation}

The proof of the Theorem \ref{Pthm2} is divided into 3 cases of evaluation $R_m(p, n):$
\begin{align*}
& Case \ 1. \ \ 5\leq m \leq \frac{n}{a(\epsilon)}, \ \ a(\epsilon) = \frac{1}{\epsilon^2}, \ \ 0<\epsilon < \frac{1}{100}, \ \ m\leq p \leq n-1 ; \\
& Case \ 2. \ \ \frac{n}{a(\epsilon)} < m \leq n - \frac{cn}{\log_2 n}, \ \ \  c\geq 7.36 , \ \ m\leq p \leq n-1; \\
& Case \ 3. \  \ n - \frac{cn}{\log_2 n} < m \leq n-1 , \ \ \  c\geq 7.36 , \ \ m\leq p \leq n-1.
\end{align*}

It was shown in the paper {\rm \cite{Odl}} that
\begin{equation}\label{eqP33}
{\mathbb P}_2(p, n) = O(p^2 2^{-n}) \ \ \mbox{as} \ \ n\to \infty, \ \ \mbox{for} \ \ 2\leq p\leq n-1;
\end{equation}

\begin{equation}\label{eqP34}
{\mathbb P}_3(p, n) =4\dbinom{p}{3}\left(\frac{3}{4}\right)^n + O\left(p^4\left(\frac{5}{8}\right)^n\right)\ \mbox{as} \ \ n\to \infty,\ \mbox{for}\ 3\leq p\leq n-1;
\end{equation}

\begin{equation}\label{eqP35}
{\mathbb P}_4(p, n) = O\left(p^4 2^{-n}\right) \ \ \mbox{as} \ \ n\to \infty, \ \ \mbox{for} \ \ 4\leq p\leq n-1.
\end{equation}

\noindent The proofs of cases 1 and 2 repeat some arguments of the papers {\rm \cite{Odl}} and {\rm \cite{IrK}}. Here we present the proofs of cases 1 and 2 for completeness of presentaton and clarification of some constants. The proof of case 3 is based on Theorem \ref{Pthm4}.

\bigskip

\subsection{Case 1 : $5\leq m \leq \frac{n}{a}, \ \  m\leq p \leq n-1 .$}

\bigskip 

\begin{lemma}\label{lmP1}
For any $\epsilon,$ $m,$ $p,$ such that $0<\epsilon < \frac{1}{100},$ $5\leq m \leq \frac{n}{a},$ where $a=a(\epsilon)=\frac{1}{\epsilon^2},$ and $m\leq p \leq n-1,$ we have
\begin{equation*}
R_m(p, n) < \left(\frac{5}{8}\right)^n(1 + \epsilon)^n \ \ \mbox{as} \ \  n \to \infty.
\end{equation*}
\end{lemma}
\proof  Let $M\in \mathcal{M}_m(m,n).$ Denote by $w_1, \ldots, w_m$ the rows of $M.$  If columns $1 \leq j_1 < \cdots < j_m \leq n$ of the matrix $M$ are linearly independent, then for each choice of $\beta_1, \ldots , \beta_m \in \{\pm 1 \},$ there will be a unique set of 
coefficients $\alpha_1, \ldots , \alpha_m$ with $j_s$th coordinate of the vector $\alpha M = ((\alpha M)_1, \ldots , (\alpha M)_n) \MYdef \alpha_1 w_1 + \cdots + \alpha_m w_m$ equals to $\beta_s,$ $s= 1, \ldots , m.$ Hence, there are at most $2^m$ sets 
$\alpha_1, \ldots, \alpha_m \in {\mathbb R}\setminus \{0\}$ such that $(\alpha M)_j = +1 \ \mbox{or} \ -1$ for $j= j_1, \ldots , j_m.$ For each fixed vector $\alpha = (\alpha_1, \ldots, \alpha_m),$ $\alpha_i \in {\mathbb R}\setminus \{0\},$ $i= 1, \ldots, m,$  probability that 
$(\alpha M)_j = +1 \ \mbox{or} \ -1$ for $j \ne j_1, \ldots , j_m,$ is at most 
$$2\cdot 2^{-m}\dbinom{m}{\left\lfloor\frac{m}{2}\right\rfloor}$$
\noindent by Erd\"os-Littlewood-Offord lemma (see {\rm \cite{Erd}}). Since all columns $ j,$ $j \ne j_1, \ldots , j_m,$ we choose independently of each other, we have
\begin{equation}\label{eqP36} 
\begin{split} 
&R_m(m,n) \leq 2^m\dbinom{n}{m}\left[2\cdot 2^{-m}\dbinom{m}{\left\lfloor\frac{m}{2}\right\rfloor}\right]^{n-m} = \\
&= 2^n\dbinom{n}{m}\left[ 2^{-m}\dbinom{m}{\left\lfloor\frac{m}{2}\right\rfloor}\right]^{n-m}.
\end{split}
\end{equation}

Taking into account (\ref{eqP29}) and (\ref{eqP36}), we get
\begin{equation}\label{eqP37}
\begin{split}
R_m(p, n) \leq & 2^n \dbinom{p}{m}\dbinom{n}{m}\left[ 2^{-m}\dbinom{m}{\left\lfloor\frac{m}{2}\right\rfloor}\right]^{n-m} .\\
\end{split}
\end{equation}

For $5\leq m \leq \frac{n}{a}, \ \  m\leq p \leq n-1,$  we have
\begin{equation*}
\begin{split}
2^n\dbinom{p}{m}\dbinom{n}{m} & \leq 2^n\dbinom{n}{m}^2 \leq 2^n\dbinom{n}{\frac{n}{a}}^2 \leq 2^n(a\cdot e)^{\frac{2n}{a}} =2^{n\left(1 + \frac{2}{a}\log_{2} a\cdot e\right)}; \\
\end{split}
\end{equation*}
\begin{equation*}
\left[ 2^{-m}\dbinom{m}{\left\lfloor\frac{m}{2}\right\rfloor}\right]^{n-m} \leq \left(\frac{5}{16}\right)^{n\left(1 - \frac{1}{a}\right)};
\end{equation*}
Thus, we have
\begin{equation*}
\begin{split}
R_m(p, n) & \leq 2^{n\left(1 + \frac{2}{a}\log_{2} a\cdot e\right)}\cdot  \left(\frac{5}{16}\right)^{n\left(1 - \frac{1}{a}\right)}  = \left(\frac{5}{8}\right)^n \left(2^{\frac{2}{a}\log_{2} a\cdot e}\left(\frac{16}{5}\right)^{\frac{1}{a}}\right)^n .
\end{split}
\end{equation*}
For any $\epsilon,$ $0<\epsilon < \frac{1}{100},$ if we take $a=a(\epsilon)=\frac{1}{\epsilon^2},$ we get
\begin{equation*}
\left(2^{\frac{2}{a}\log_{2} a\cdot e}\left(\frac{16}{5}\right)^{\frac{1}{a}}\right) < 1+ \epsilon.
\end{equation*}
Hence, for any $\epsilon,$ $m,$ $p,$ such that $0<\epsilon < \frac{1}{100},$ $5\leq m \leq \epsilon^2 n, \ \  m\leq p \leq n-1,$ we have
\begin{equation}\label{eqP38}
R_m(p, n) < \left(\frac{5}{8}\right)^n (1+ \epsilon)^n \ \ \mbox{as} \ \ n \to \infty.
\end{equation}

\begin{flushright} {\sc Q.E.D.} \end{flushright}

\bigskip

\subsection{Case 2 : $ \epsilon^2 n <m  \leq n - \frac{cn}{\log_2n}, \   0< \epsilon < \frac{1}{100}, \   c  \geq 7.36, \   m\leq p \leq n-1 .$}

\bigskip 

\begin{lemma}\label{lmP2}
For any $\epsilon,$ $m,$ $p,$ such that $0<\epsilon < \frac{1}{100},$ $\epsilon^2 n < m  \leq n - \frac{cn}{\log_2n}, $ where $c \geq 7.36,$ and $m\leq p \leq n-1,$ we have
\begin{equation*}
R_m(p, n) = o\left(\left(\frac{5}{8}\right)^n\right) \ \ \mbox{as} \ \  n \to \infty.
\end{equation*}
\end{lemma}
\proof  Using arguments from the first case, we have:
\begin{equation*} 
\begin{split} 
R_m(m,n) &\leq 2^m\dbinom{n}{m}\left[2\cdot 2^{-m}\dbinom{m}{\left\lfloor\frac{m}{2}\right\rfloor}\right]^{n-m} \leq 2^{2n}\left[ 2^{-m}\dbinom{m}{\left\lfloor\frac{m}{2}\right\rfloor}\right]^{n-m} \leq \\
&\leq 2^{2n}\left[\sqrt{\frac{2}{\pi \epsilon^2}}n^{-\frac{1}{2}}\right]^{\frac{cn}{\log_2 n}} = 2^{2n - \frac{cn}{2}}\left[\left(\frac{2}{\pi \epsilon^2}\right)^{\frac{c}{2 \log_{_2} n}}\right]^n.
\end{split}
\end{equation*}
Then,
\begin{equation}\label{eqP39}
\begin{split}
R_m(p,m) &\leq \dbinom{p}{m}R_m(m,n) \leq \dbinom{p}{m}2^{2n - \frac{cn}{2}}\left[\left(\frac{2}{\pi \epsilon^2}\right)^{\frac{c}{2 \log_2 n}}\right]^n \leq \\
&\leq 2^{3n - \frac{cn}{2}}\left[\left(\frac{2}{\pi \epsilon^2}\right)^{\frac{c}{2 \log_2 n}}\right]^n =  o\left(\left(\frac{5}{8}\right)^n\right) \ \ \mbox{for $c \geq 7,36$.}
\end{split}
\end{equation}

\begin{flushright} {\sc Q.E.D.} \end{flushright}

\bigskip

\subsection{Case 3 : $  \ n - \frac{cn}{\log_2 n} < m \leq n-1 , \ \ \  c\geq 7.36 , \ \ m\leq p \leq n-1.$}

\bigskip

\begin{lemma}\label{lmP3}
For any $m,$  $n - \frac{cn}{\log_{2}n} < m \leq n-1,$ where $c \geq 7.36,$ and $p,$ $m\leq p \leq n-1,$ we have
\begin{equation*}
R_m(p, n) = \left(\frac{1}{2} + o_n(1)\right)^n \ \ \mbox{as} \ \  n \to \infty.
\end{equation*}
\end{lemma}
\proof  Let $M\in \mathcal{M}_m(m,n)$ and $M(j_1, \ldots ,j_{m+1})$ be its $m \times (m+1)$-submatrix with columns $j_1 < \ldots < j_{m+1}.$ Denote by $\mathcal{M}_m(m,n; j_1, \ldots , j_{m+1})$ the set
\begin{equation*}
\begin{split}
&\mathcal{M}_m(m,n; j_1, \ldots , j_{m+1}) \MYdef \\
&\MYdef \{M\in \mathcal{M}_m(m,n)\mid  M(j_1, \ldots ,j_{m+1}) \in  \mathcal{M}_m(m,m+1) \}.
\end{split}
\end{equation*}
Then
\begin{equation}\label{eqP40}
\mathcal{M}_m(m,n) \subset \bigcup_{1\leq j_1 < \ldots < j_{m+1}\leq n} \mathcal{M}_m(m,n; j_1, \ldots , j_{m+1}) .
\end{equation}
On the other hand, by Theorem \ref{Pthm4} we have:
\begin{equation}\label{eqP41}
\begin{split}
R_m(m, m+1) &= \frac{|\mathcal{M}_m(m,m+1)|}{2^{m(m+1)}} =\\
&=\frac{|KSO(m, m+1) \cap \langle E_m \rangle^{\times m}_{\ne 0}|}{2^{m^2}} \leq \frac{m^2}{2^{m-1}} .
\end{split}
\end{equation}
From (\ref{eqP40}), (\ref{eqP41}), and (\ref{eqP29}), we get
\begin{equation*}
\begin{split}
R_m(p, n) &\leq \dbinom{p}{m}\dbinom{n}{m+1}R_m(m, m+1) \leq \dbinom{n}{m}^2 \frac{m^2}{2^{m-1}} \leq \\
&\leq \dbinom{n}{\frac{cn}{\log_2 n}}^2\frac{m^2}{2^{m-1}}\leq \left(\frac{e\log_2 n}{c}\right)^{\frac{2cn}{\log_2 n}}\cdot \frac{n^2}{2^{n-\frac{cn}{\log_2 n}}} = \left(\frac{1}{2}+ o_n(1)\right)^n.
\end{split}
\end{equation*}
\begin{flushright} {\sc Q.E.D.} \end{flushright}

Now Theorem \ref{Pthm2} follows from (\ref{eqP31}), (\ref{eqP32}),  (\ref{eqP33}), (\ref{eqP34}), (\ref{eqP35}), Lemma \ref{lmP1},  Lemma \ref{lmP2}, and  Lemma~\ref{lmP3}.

\begin{flushright} {\sc Q.E.D.} \end{flushright}

\bigskip

\noindent Faculty of Mechanics and Mathematics, Lomonosov Moscow State University, Moscow, Russian Federation

\bigskip

{\it E-mail address}:\ \ \ anwar.irmatov@math.msu.ru


\begin{thebibliography}{99}



\bibitem{BVW}

J.~Bourgain, V.~H.~Vu, P.~M.~Wood, On the singularity probability of discrete randon matrices, {\it Journal of Functional Analysis}, 258~(2010), 559--663.  DOI: https://doi.org/10.1016/j.jfa.2009.04.016



\bibitem{Erd}

P.~Erd\" os, On a lemma of Littlewood and Offord, {\it Bull. Amer. Math. Soc.} 54~(1945), 898--902. MR0014608

\bibitem{Ir1}

A.~A.~Irmatov, On the number of threshold functions, {\it Discrete Math. Appl.}, 3(4)~(1993), 429--432. DOI: https://doi.org/10.1515/dma-1993-0407



\bibitem{IrK}

A.~A.~Irmatov, \v Z.~D.~Kovijani\' c, On the asymptotics of logarithm of the number of threshold K-logic functions, {\it Discrete Math. Appl.}, 8(4)~(1998), 331--355. DOI: https://doi.org/10.1515/dma.1998.8.4.331

\bibitem{Ir4}

A.~A.~Irmatov, Singularity of $\{\pm 1\}$-matrices and asymptotics of the number of threshold functions, https://arxiv.org/abs/2004.03400  DOI: https://doi.org/10.48550/arXiv.2004.03400

\bibitem{Ir5}

A.~A.~Irmatov, Asymptotics of the number of threshold functions and the singularity probability of random $\{\pm 1\}$-matrices, Doklady Rossijskoj Akademii Nauk. Mathematika, Informatika, Processy Upravlenia, 2020, Volume 492, Pages 89--91
DOI: https://doi.org/10.31857/S2686954320030091



\bibitem{KKS}

J.~Kahn, J.~Koml\' os, and E.~Szemer\' edi, On the probability that a random $\pm 1$-matrix is singular, {\it Journal of the Amer. Math.Soc.}, Vol.~8, Num.~1, Jan. 1995, 223--240. MR1260107, DOI: https://doi.org/10.1090/S0894-0347-1995-1260107-2

\bibitem{KaSo}

I.~Kanter and H.~Sompolinsky, Associative recall of memory without errors, Phys Rev A Gen Phys. 1987 Jan 1;35(1):380-392. doi: 10.1103/physreva.35.380

\bibitem{Ko1}

J.~Koml\' os, On the determinant of (0,1) matrices, {\it Studia Sci. Math. Hungar.}, 2~(1967), 7--21. MR0221962

\bibitem{Ko2}

J.~Koml\' os, Manuscript (1977). In: {\it B.Bollob\' as (Ed.) ``Random Graphs'', Academic Press, New York/London}, 1985, 347--350. 


\bibitem{Odl}

A.~M.~Odlyzko, On subspaces spanned by random selection of $\pm 1$ vectors, {\it J. Combin. Theory Ser. A}, 47~(1988), 124--133. DOI: https://doi.org/10.1016/0097-3165(88)90046-5

\bibitem{Rot}

G.-C.~Rota, On the Foundations of Combinatorial Theory I. Theory of M\"obius Functions, Z. Wahrscheinlichkeitstheorie 2 (1964), 340-368.


\bibitem{Sch}

L.~Schl\" afli, Gesammelte Mathematische Abhandlungen I, {\it Verlag Birkh\" auser, Springer Basel AG}, 1950, 209--212.

\bibitem{TV2}

T.~Tao and V.~Vu, On the singularity of random Bernoulli matrices, {\it Journal of Amer. Math. Soc.}, Vol.20, Num.3, July 2007, 603--628.  MR2291914, DOI: https://doi.org/10.1090/S0894-0347-07-00555-3

\bibitem{Tik}

K.~Tikhomirov, Singularity of random Bernoulli matrices, {\it Annals of Math.}, Vol.191~(2020), Issue 2, 593--634 DOI: https://doi.org/10.4007/annals.2020.191.2.6


\bibitem{Zas}

T.~Zaslavsky, Facing up to Arrangements: Face-count Formulas for Partitions of Space by Hyperplanes, {\it Mem. Amer. Math. Soc.}, 154, Amer.Math.Soc., Providence, RI, 1975. MR 50:9603

\bibitem{Zue}

Yu.~A.~Zuev, Asymptotics of the logarithm of the number of threshold functions of the algebra of logic, {\it Sov. Math., Dokl..}, 39~(1989), 512--513. Zbl 0693.94010

 

\end{thebibliography}
\end{document}